\documentclass{amsart}

\usepackage[latin1]{inputenc}
\usepackage{latexsym}
\usepackage{amsfonts}
\usepackage{amssymb}
\usepackage{color}
\usepackage{dsfont}
\usepackage{varioref}
\usepackage{epsfig}
\usepackage{hyperref}

\newtheorem{Lem}{Lemma}

\newtheorem{Thm}{Theorem}

\newtheorem{Cor}[Thm]{Corollary}

\newcommand{\dN}{{\mathds N}}
\newcommand{\dP}{{\mathds P}}
\newcommand{\dC}{{\mathds C}}
\newcommand{\dR}{{\mathds R}}

\newcommand{\smcdot}{{\textup{$\cdot$}}}

\newcommand{\dontShowThis}[1]{}

\begin{document}

\title[Real Line Arrangements and Surfaces with Many Real
  Nodes]{Real Line Arrangements and Surfaces\\with Many Real
  Nodes%\\--- PRELIMINARY VERSION ---
} 
\begin{abstract}
A long standing question is if maximum number $\mu(d)$ of nodes 
on a surface of degree $d$ in $\dP^3(\dC)$ can be achieved by
a surface defined over the reals which has only real singularities.  
The currently best known asymptotic lower bound, $\mu(d) \gtrapprox 
\frac{5}{12}d^3$, is provided by
Chmutov's construction from 1992 which gives surfaces whose nodes have  
non-real coordinates.  

Using explicit constructions of certain real line arrangements we show 
that Chmutov's construction can be adapted to give only real
singularities.   
All currently best known constructions which
exceed Chmutov's lower bound (i.e., for $d=3,4,\cdots,8,10,12$)
can also be realized with only real singularities.  
Thus, our result shows that, up to now, all known lower bounds for $\mu(d)$
can be attained with only real singularities. 

We conclude with an application of the theory of real line
arrangements which shows that our arrangements are aymptotically the best
possible ones.
This proves a special case of a conjecture of Chmutov. 
\end{abstract}
\date{\today}
\author{Sonja Breske}
%\address{Johannes Gutenberg Universit\"at Mainz, Germany}
%\email{???} 
\author{Oliver Labs} 
%\address{Johannes Gutenberg Universit\"at Mainz, Germany}
%\email{Labs@Mathematik.Uni-Mainz.de, }
\author{Duco van Straten} 
\address{Institut f\"ur Mathematik, Johannes Gutenberg Universit\"at Mainz,
  D-55099 Mainz, Germany}
\email{\{Breske,Labs,Straten\}@Mathematik.Uni-Mainz.de}

%    General info

\maketitle

%\begin{center}
%Affiliation: 
%Johannes Gutenberg Universit\"at Mainz, Institut f\"ur Mathematik,\\
%D-55099 Mainz, Germany, \\
%E-Mail: {\{Breske,Labs,Straten\}@Mathematik.Uni-Mainz.de}\\
%Corresponding author: O.\ Labs, Tel: +49 6131 3923335, Fax: +49 6131
%3920915. 
%\end{center}

%%%%%%%%%%%%%%%%%%%%%%%%%%%%%%%%%%%%%%%%%%%%%%%%%%%%%%%%%%%%%%%%%%
%
%
%
%%%%%%%%%%%%%%%%%%%%%%%%%%%%%%%%%%%%%%%%%%%%%%%%%%%%%%%%%%%%%%%%%%
\section*{Introduction}

A node (or $A_1$ singularity) in $\dC^3$ is a singular point which can be
written in the form $x^2+y^2+z^2=0$ in some local coordinates. 
We denote by $\mu(d)$ the maximum possible number of nodes 
on a surface in $\dP^3(\dC)$. 
The question of determining $\mu(d)$ has a long and rich history. 
Currently, $\mu(d)$ is only known for $d=1,2,\dots,6$ (see \cite{bar65,
  jafrub66} for sextics and \cite{labs99} for a recent improvement
for septics). 

In this paper, we consider the relationship between $\mu(d)$ and the maximum
possible number of real nodes on a surface 
in $\dP^3(\dR)$ which we denote by $\mu^\dR(d)$. 
Obviously, $\mu^\dR(d) \le \mu(d)$, but do we even have $\mu^\dR(d) = \mu(d)$?
In other words: Can the maximum number of nodes be achieved with
real surfaces with real singularities?

The previous question arises naturally because all results in low degree $d\le12$
  suggest that it could be true (see \cite{bar65, labs99, endrOct,
  endrHistory, ales12} and table \ref{tabBounds}).   
But the best known asymptotic lower bound, $\mu(d)\gtrapprox \frac{5}{12}d^3$, 
follows from Chmutov's construction \cite{chmuP3} which yields only
singularities with non-real coordinates. 
In this paper, we show that his construction can be adapted to give surfaces
with only real singularities (see table \vref{tabBounds}). 
In the real case we can distinguish between two types of nodes,
\emph{conical nodes} ($x^2+y^2-z^2=0$) and \emph{solitary points}
($x^2+y^2+z^2=0$): Our construction produces only conical nodes.

Notice that in general there are no better real upper bounds for $\mu^\dR(d)$
known than the well-known complex ones of Miyaoka \cite{miyp3} and Varchenko
\cite{varBound}. 
But for solitary points there exist better bounds via the relation to the
zero$^\textup{th}$ Betti number (see e.g., \cite{kharTopRealMan}).  
%For solitary points there are better bounds because of its relation to the
%zero$^\textup{th}$ Betti number:
%E.g., a real sextic curve in the real plane
%cannot have nine real cusps which is the complex maximum number. 
E.g., Rohn 
%(or maybe already R.W.H.T.\ Hudson) 
showed in 1913 that a real
quartic surface in $\dP^3(\dR)$ cannot 
have more than $10$ solitary points although it can have $16$ conical nodes.   
We show a real upper bound of $\approx \frac56 d^2$ for the
maximum number of critical points on two levels of real simple line
arrangements consisting of $d$ lines.  
In \cite{chmuCritVals}, Chmutov conjectured this to be the maximum number
for all complex plane curves of degree $d$. 
He also noticed \cite{chmuP3} that such a bound directly implies an upper bound for the
number of real nodes of certain surfaces.
Our upper bound shows that our examples are asymptotically the best possible
real line arrangements for this purpose.
%The proof uses the relation to two-colorings of such arrangements.

\begin{table}[htbp]
%  \hspace*{-0.75cm}
  \begin{tabular}{|r|@{\;}c@{\;}|@{\;}c@{\;}|@{\;}c@{\;}|@{\;}c@{\;}|@{\;}c@{\;}|@{\;}c@{\;}|@{\;}c@{\;}|@{\;}c@{\;}|@{\;}c@{\;}|@{\;}c@{\;}|@{\;}c@{\;}|@{\;}c@{\;}|@{\;}c@{\;}|@{\;}c@{\;}|}
    \hline 
    \rule{0pt}{1.2em}$d$ & $1$ & $2$ & $3$ & $4$ & $5$ & $6$ & $7$ & $8$ & $9$ &
    $10$ & $11$ & $12$ & $13$ & $d$\\[0.1em]
    \hline 
    \rule{0pt}{1.2em}$\mu(d), \mu^\dR(d)\le$ & $0$ & $1$ & $4$ & $16$ & $31$ & $65$ & $104$ & $174$ &
    $246$ & $360$ & $480$ & $645$ & $832$ & $\frac{4}{9}d(d-1)^2$\\[0.1em]
    \hline 
    \rule{0pt}{1.2em}$\mu(d), \mu^\dR(d)\ge$ & $0$ & $1$ & $4$ & $16$ & $31$ & $65$ & $99$ & $168$ & $\mathbf{216}$ &
    $345$ & $\mathbf{425}$ & $600$ & $\mathbf{732}$ & $\approx \mathbf{\frac{5}{12}d^3}$ \\[0.1em]
    \hline 
  \end{tabular}\\[0.3em]
  \caption{The currently known bounds for the maximum number $\mu(d)$ (resp.\
    $\mu^\dR(d)$) of nodes on a surface of degree $d$ in $\dP^3(\dC)$
    (resp.\ $\dP^3(\dR)$) are equal. 
    The bold numbers indicate in which cases our result improves the previously known
    lower bound for $\mu^\dR(d)$.}
  \label{tabBounds}
\end{table}

%%%%%%%%%%%%%%%%%%%%%%%%%%%%%%%%%%%%%%%%%%%%%%%%%%%%%%%%%%%%%%%%%%%%%%
%
%
%%%%%%%%%%%%%%%%%%%%%%%%%%%%%%%%%%%%%%%%%%%%%%%%%%%%%%%%%%%%%%%%%%
\section{Variants of Chmutov's Surfaces with Many Real Nodes}

\begin{figure}[htbp]
  \begin{center}
  \begin{tabular}{ccc}
    \includegraphics[width=1.4in]{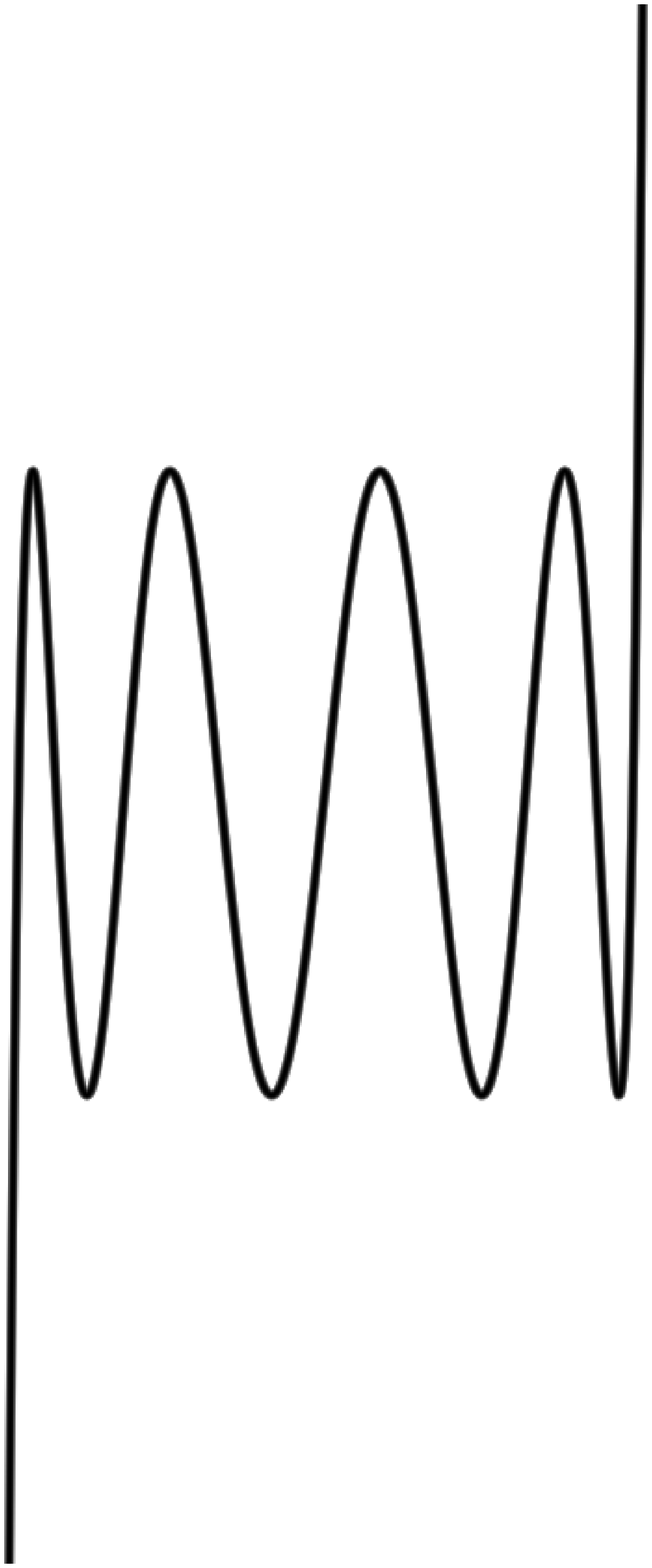}
    & 
    \includegraphics[width=1.4in]{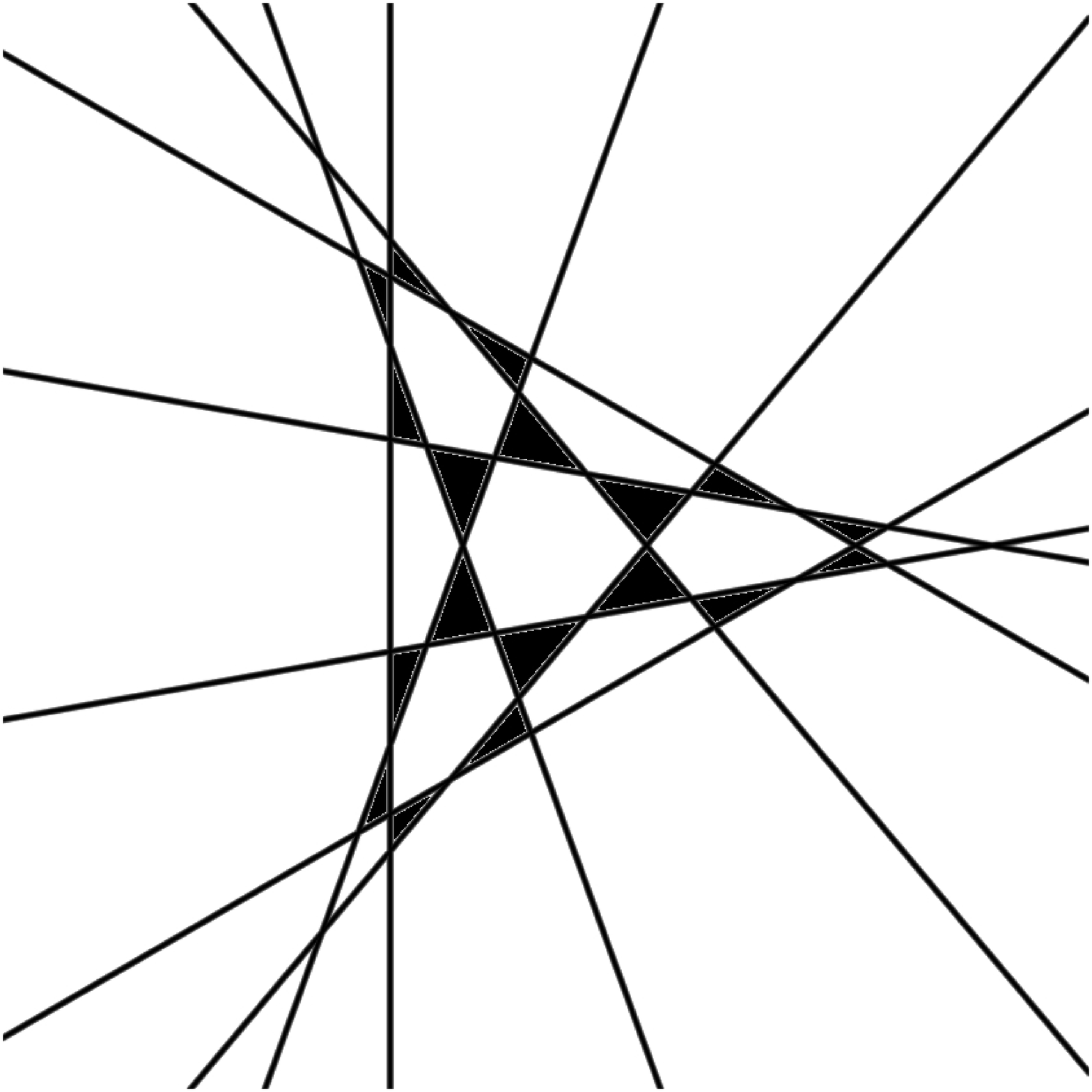}
    &
    \includegraphics[width=1.4in]{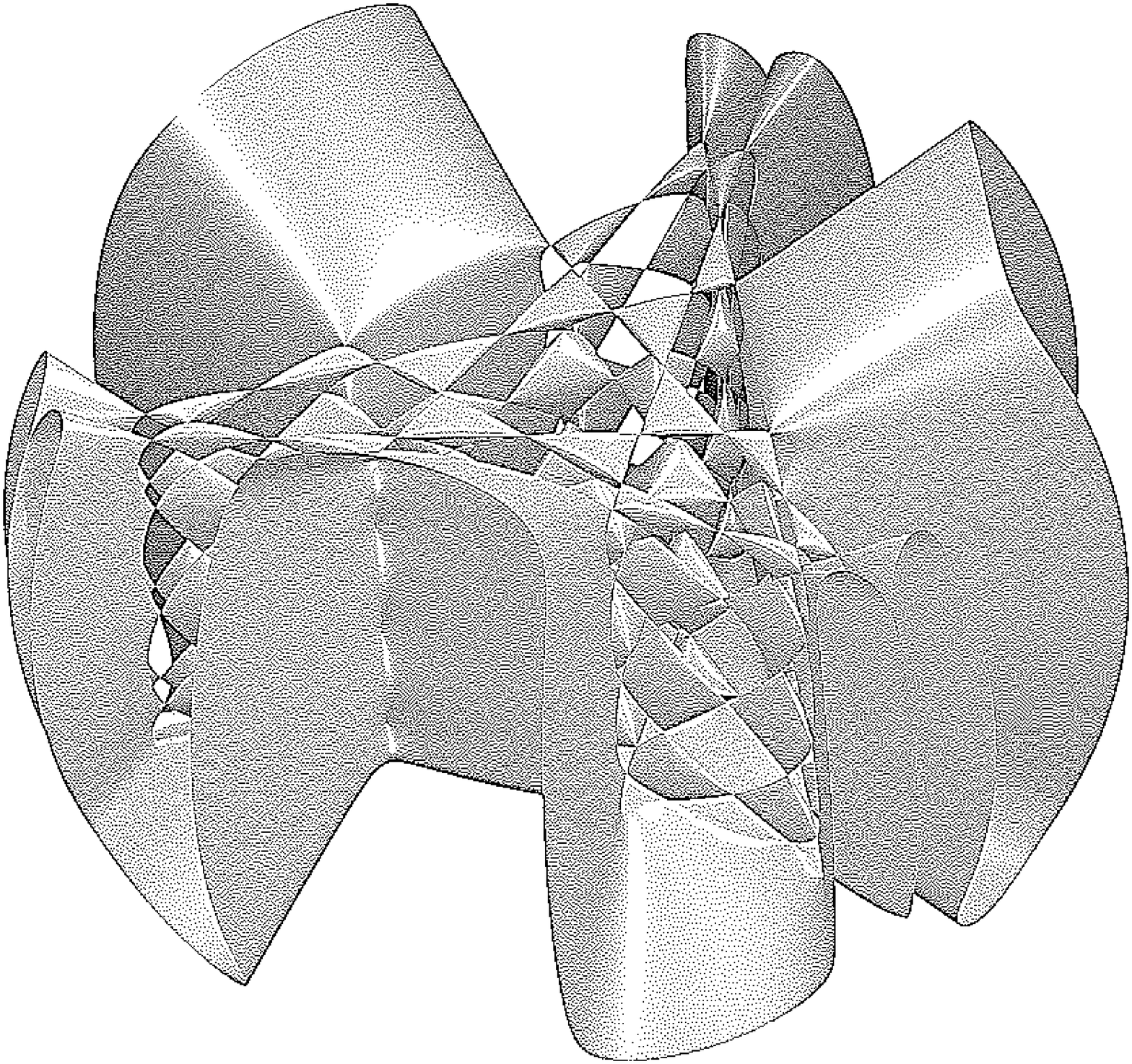}
  \end{tabular}
  \caption{For degree $d=9$ we show the Tchebychev polynomial $T_9(z)$, the
    real folding polynomial $F^{A_2}_{\dR,9}(x,y)$ associated to
    the root system $A_2$, and
    the surface $\textup{Chm}^{A_2}_{\dR,9}(x,y,z)$.
    The bounded regions in which $F^{A_2}_{\dR,9}(x,y)$ takes negative values are marked
    in black.} 
  \end{center}
  \label{figTchebFold}
\end{figure}

Let $T_d(z)\in\dR[z]$ be the Tchebychev polynomial of degree $d$ with critical
values $-1$ and $+1$ (see fig.\ \vref{figTchebFold}), i.e.:
$T_0(z) := 1$, \ $T_1(z) := z$, \ $T_d(z) := 2\smcdot z\smcdot T_{d-1}(z) -
T_{d-2}(z)$ for $d\ge2$.
Chmutov \cite{chmuP3} uses them together with the so-called \emph{folding
polynomials} $F^{A_2}_{d}(x,y)\in\dR[x,y]$ associated to the root-system $A_2$ to
construct surfaces $\textup{Chm}_{d}^{A_2}(x,y,z) := F^{A_2}_{d}(x,y) +
\frac{1}{2}(T_d(z) + 1)$ with many nodes.
These folding polynomials are defined as follows:
{\small 
  \begin{equation}\label{eqnFA2d}
    F^{A_2}_d(x,y) := 2+\det\left(
      \begin{array}{c@{\,}c@{\,}c@{\,}c@{\,}c@{\,}c@{\,}c}
        \rule{0pt}{1.4em}x & 1 & 0 & \cdots & \cdots & \cdots & 0\\[-0.6em]
        \rule{0pt}{1.4em}2y & x & \ddots & \ddots &  &  & \vdots \\[-0.6em]
        \rule{0pt}{1.4em}3 & y & \ddots & \ddots & \ddots & & \vdots \\[-0.6em]
        \rule{0pt}{1.4em}0 & 1 & \ddots & \ddots & \ddots & \ddots & \vdots \\[-0.6em]
        \rule{0pt}{1.4em}\vdots & \ddots & \ddots & \ddots & \ddots & \ddots & 0\\[-0.6em]
        \rule{0pt}{1.4em}\vdots &  & \ddots & \ddots & \ddots & \ddots & 1\\[-0.6em]
        \rule{0pt}{1.4em}0 & \cdots & \cdots & 0 & 1 & y & x
      \end{array}
    \right) + 
    \det\left(
      \begin{array}{c@{\,}c@{\,}c@{\,}c@{\,}c@{\,}c@{\,}c}
        \rule{0pt}{1.4em}y & 1 & 0 & \cdots & \cdots & \cdots & 0\\[-0.6em]
        \rule{0pt}{1.4em}2x & y & \ddots & \ddots &  &  & \vdots \\[-0.6em]
        \rule{0pt}{1.4em}3 & x & \ddots & \ddots & \ddots & & \vdots \\[-0.6em]
        \rule{0pt}{1.4em}0 & 1 & \ddots & \ddots & \ddots & \ddots & \vdots \\[-0.6em]
        \rule{0pt}{1.4em}\vdots & \ddots & \ddots & \ddots & \ddots & \ddots & 0\\[-0.6em]
        \rule{0pt}{1.4em}\vdots &  & \ddots & \ddots & \ddots & \ddots & 1\\[-0.6em]
        \rule{0pt}{1.4em}0 & \cdots & \cdots & 0 & 1 & x & y
      \end{array}
    \right).\end{equation}
}

%\begin{equation}\label{eqnFA2d}
%  F^{A_2}_{d}(x,y):=P^{A_2}_{1,d}(x,y)+P^{A_2}_{2,d}(x,y)+2,
%\end{equation}
%where the $P^{A_2}_{i,d}$ are recursively defined as follows (see \cite[p.\
%66--67]{sonjaChmuVars}):
%\begin{center}\begin{tabular}{l@{\qquad}l@{\qquad}l}
%$P^{A_2}_{1,1}(x,y):=x$, & $P^{A_2}_{1,2}(x,y):=x^2-2y$, & $P^{A_2}_{1,3}(x,y):=x^3-3xy+3$, \\[0.2em]
%$P^{A_2}_{2,1}(x,y):=y$, & $P^{A_2}_{2,2}(x,y):=y^2-2x$, & $P^{A_2}_{2,3}(x,y):=y^3-3xy+3$,
%\end{tabular}\end{center}
%\begin{align*}
%P^{A_2}_{1,d}(x,y):=x\smcdot P^{A_2}_{1,(d-1)}(x,y)-y\smcdot P^{A_2}_{1,(d-2)}(x,y)+P^{A_2}_{1,(d-3)}(x,y),\\
%P^{A_2}_{2,d}(x,y):=y\smcdot P^{A_2}_{2,(d-1)}(x,y)-x\smcdot P^{A_2}_{2,(d-2)}(x,y)+P^{A_2}_{2,(d-3)}(x,y).
%\end{align*}  
The $F^{A_2}_{d}(x,y)$ have critical points with only three different critical
values: $0$, $-1$, and $8$.  
Thus, the surface $\textup{Chm}_{d}^{A_2}(x,y,z)$ is singular exactly at those
points at which the critical values of $F^{A_2}_{d}(x,y)$ and
$\frac{1}{2}(T_d(z) + 1)$ sum up to zero (i.e., either both are $0$ or the
first is $-1$ and the second is $+1$). 

Notice that the plane curve defined by $F^{A_2}_{d}(x,y)$ consists in fact of
$d$ lines. 
But these are not real lines and the critical points of this folding
polynomial also have non-real coordinates.
It is natural to ask whether there is a real line arrangement which leads to 
the same number of critical points. 
The term \emph{folding polynomials} was introduced in \cite{witFoldPoly}
(here we use a slightly different definition).  
In his article, Withers also described many of their properties, 
but it was Chmutov \cite{chmuP3} who noticed that $F^{A_2}_{d}(x,y)$ has only
few different critical values. 
In \cite{sonjaChmuVars}, the first author computed the critical points of the
other folding polynomials. 
Among these, there are the following examples which are the real line
arrangements we have been looking for (see \cite[p.\ 87--89]{sonjaChmuVars}): 

We define the real folding polynomial $F^{A_2}_{\dR,d}(x,y)\in\dR[x,y]$
associated to the root system $A_2$ as (see also fig.\ \vref{figTchebFold}) 
\begin{equation}
 F^{A_2}_{\dR,d}(x,y) := F^{A_2}_{d}(x+iy, \ x-iy),
\end{equation}
where $i$ is the imaginary
number.  
It is easy to see that the $F^{A_2}_{\dR,d}(x,y)$ have indeed real
coefficients. 
The numbers of critical points are the same as those of $F^{A_2}_{d}(x,y)$;
but now they have real coordinates as the following lemma shows:

%\begin{Def}
%We define the real folding polynomial $F^{A_2}_{\dR,d}(x,y)$ of degree $d\in\dN$
%associated to the root system $A_2$ by
%$F^{A_2}_{\dR,d}(x,y):=P^{A_2}_{1,d}+\overline{P^{A_2}_{1,d}}+2$, where
%$\overline{\rule{0pt}{0.8em}\ \cdot \ }$ denotes the complex conjugation and
%\marginpar{Sonja: simplify this recursive definition by deriving it from the
%  other one (without $i$)??? bitte auch angeben, auf welchen Seiten Deiner
%  Arbeit die entsprechenden rekursiven Definitionen zu finden sind!!!} 
%\begin{align*}
%  &P^{A_2}_{1,1}:=x+i\smcdot y, \quad P^{A_2}_{1,2}:= x^2-y^2-2x+i\smcdot(2xy+2y), \\   
%  &P^{A_2}_{1,3}:=x^3-3xy^2-3x^2-3y^2+3+i\smcdot(3x^2y-y^3).
%\end{align*}
%When writing $P^{A_2}_{1,d}$ in the form $P^{A_2}_{1,d}=a_d+i\smcdot b_d$ then the
%$d^\textup{th}$ polynomial is given by:
%\begin{align*}
%  &P^{A_2}_{1,d}:=(x+i\smcdot y)(a_{d-1}-a_{d-2})+(-y+i\smcdot x)(b_{d-1}+b_{d-2})+a_{d-3}+b_{d-3}.  
%\end{align*}
%\end{Def}

\begin{Lem}
  The real folding polynomial $F^{A_2}_{\dR,d}(x,y)$ associated to the root 
  system $A_2$ has $d\choose2$ real critical points with critical value 
  $0$ and 
  \begin{align}
      {1\over3}d^2 - d & \quad \textup{if} \ d\equiv 0 \mod 3, & 
      {1\over3}d^2-d+{2\over3} & \quad  \textup{otherwise}  
    \end{align}
  real critical points with critical value $-1$.
  The other critical points also have real coordinates and have critical value
  $8$. 
\end{Lem}
\begin{proof}
We proceed similar to the case discussed in \cite{chmuP3}, see
\cite[p.\ 87--95]{sonjaChmuVars} for details.  
To calculate the critical points of the real folding polynomial $F^{A_2}_{\dR,d}$, we
use the map $h^1: \dR^2\to\dR^2$, defined by 
%{\small
\[(u,v) \mapsto \bigl(\cos(2\pi(u+v))+\cos(2\pi u)+\cos(2\pi v), \ \ 
\sin(2\pi(u+v))-\sin(2\pi
  u)-\sin(2\pi v)\bigr).\]
%}
This is in fact just the real and imaginary part of the first component of the
generalized cosine $h$ considered by Withers \cite{witFoldPoly} and Chmutov
\cite{chmuP3}.  
It is easy to see that $h^1$ is a coordinate
change if $u-v>0,$ $u+2v>0$, and $2u+v<1$.
It transforms the polynomial $F^{A_2}_{\dR,d}$ into the function $G^{A_2}_{d}:
\ \dR^2\to\dR^2$, defined by 
\[
 G^{A_2}_{d}(u,v) := F^{A_2}_{\dR,d}(h^1(u,v)) = 2\cos(2\pi du)+2\cos(2\pi 
 dv)+2\cos(2\pi d(u+v))+2.  
\] 
The calculation of the critical points of $G^{A_2}_{d}$ is exactly the same as
the one performed in \cite{chmuP3}. As the function $G^{A_2}_{d}$ has
$(d-1)^2$ distinct real critical points in 
the region defined by $u-v>0,$ $u+2v>0$, and $2u+v<1$, the images of these
points under the map $h^1$ are all the critical points of the real folding
polynomial $F^{A_2}_{\dR,d}$ of degree $d$.  
In contrast to \cite{chmuP3}, we get real critical points because
$h^1$ is a map from $\dR^2$ into itself. 
%We proceed similar to the case discussed in \cite{chmuP3}, see
%\cite[p.\ 87--95]{sonjaChmuVars} for details.  
%To calculate the critical points of the polynomial $F^{A_2}_{\dR,d}$, we
%consider the map 
%\[h^1: \ \dR^2\to\dC, \ 
%(u,v) \mapsto \left(e^{-2\pi iu}+e^{-2\pi iv}+e^{2\pi i(u+v)}\right) 
%\]
%which is the first component of the generalized cosine used by Withers
%\cite{witFoldPoly}. 
%Then $\widetilde{h^1}(u,v) :=
%\bigl(\mathfrak{Re}(h^1(u,v)),\mathfrak{Im}(h^1(u,v))\bigr)$ is a map from
%the real plane into itself. 
%It is easy to see that $\widetilde{h^1}$ is a coordinate
%change if $u-v>0,$ $u+2v>0$ and $2u+v<1$.
%It transforms the polynomial $F^{A_2}_{\dR,d}$ into the function $G^{A_2}_{d}:
%\ \dR^2\to\dR^2$, defined by 
%\[
% G^{A_2}_{d}(u,v) := F^{A_2}_{\dR,d}(\widetilde{h^1}(u,v)) = 2\cos(2\pi du)+2\cos(2\pi 
% dv)+2\cos(2\pi d(u+v))+2.  
%\] 
%The calculation of the critical points of $G^{A_2}_{d}$ is exactly the same as
%the one performed in \cite{chmuP3}. As the function $G^{A_2}_{d}$ has $(d-1)^2$ distinct real critical points in
%the region defined by $u-v>0,$ $u+2v>0$ and $2u+v<1$, the images of these
%points under the map $\widetilde{h^1}$ are all the critical points of the real folding
%polynomial $F^{A_2}_{\dR,d}$ of degree $d$.  
%In contrast to \cite{chmuP3}, we get real critical points because
%$\widetilde{h^1}$ is a map from $\dR^2$ into itself.  
\end{proof}

None of the other root systems yield more critical points on two levels. 
But as mentioned in \cite{labsVarChmu}, the real folding polynomials
associated to the root system $B_2$ give hypersurfaces in $\dP^n$, $n\ge5$,
which improve the previously known lower bounds for the maximum number of
nodes in higher dimensions slightly (see \cite{labsVarChmu};
\cite{sonjaChmuVars} gives a detailed discussion of all these folding
polynomials and their critical points).   

The lemma immediately gives the following variant of Chmutov's nodal
surfaces:  
\begin{Thm}\label{thmChmRd}
  Let $d\in\dN$. 
  The real projective surface of degree $d$ defined by
  \begin{equation}
    \textup{Chm}_{\dR,d}^{A_2}(x,y,z) \ := \ F^{A_2}_{\dR,d}(x,y) +
    \frac{1}{2}(T_d(z) + 1)\ \in \ \dR[x,y,z]
  \end{equation}
  has the following number of real nodes:
  \begin{equation}\label{eqnChmuSurfNoNo}
    \begin{array}{ll}
      \frac{1}{12}\left(5d^3-13d^2+12d\right) & \textup{if} \ d\equiv 0 \mod 6, \\[0.4em]
      \frac{1}{12}\left(5d^3-13d^2+16d-8\right) & \textup{if} \ d\equiv 2,4 \mod 6, \\[0.4em]
      \frac{1}{12}\left(5d^3-14d^2+13d-4\right) & \textup{if} \ d\equiv 1,5 \mod 6, \\[0.4em]
      \frac{1}{12}\left(5d^3-14d^2+9d\right) & \textup{if} \ d\equiv 3 \mod 6.
    \end{array}
  \end{equation}
\end{Thm}

These numbers are the same as the numbers of complex nodes of Chmutov's 
surfaces $\textup{Chm}_{d}^{A_2}(x,y,z)$. 
To our knowledge, the result gives new lower bounds for the maximum number
$\mu^\dR(d)$ of real singularities on a surface of degree $d$ in $\dP^3(\dR)$
for $d=9,11$ and $d\ge13$, see table \vref{tabBounds}.  
Notice that all best known lower bounds for $\mu^\dR(d)$ are attained by
surfaces with only conical nodes which is not astonishing in view of the upper
bounds for solitary points mentioned in the introduction.

%%%%%%%%%%%%%%%%%%%%%%%%%%%%%%%%%%%%%%%%%%%%%%%%%%%%%%%%%%%%%%%%%%%%%%
%
%
%
%%%%%%%%%%%%%%%%%%%%%%%%%%%%%%%%%%%%%%%%%%%%%%%%%%%%%%%%%%%%%%%%%%
\section{On Two-Colorings of Real Simple Line Arrangements}

The real folding polynomials $F_{\dR,d}^{A_2}(x,y)$ used in the previous
section are in fact \emph{real simple (straight) line arrangements} in $\dR^2$,
i.e.,  lines no three of which meet in a point. 
%This already follows from the fact that $d\choose2$ of its critical points
%have critical value $0$. 
Such arrangements can be $2$-colored in a natural way (see fig.\
\vref{figTchebFold}):
We label in black those regions (\emph{cells}) of $\dR^2 \ \backslash \ 
\{F_{\dR,d}^{A_2}(x,y)=0\}$ in which $F_{\dR,d}^{A_2}(x,y)$ 
takes negative values, the others in white. 
The bounded black regions in fig.\ \ref{figTchebFold} contain exactly one
critical point with critical value $-1$ each. 

Harborth has shown in \cite{harborthTwoCol} that the maximum number $M_b(d)$
of black cells in such real simple line arrangements of $d$ lines satisfies:
\begin{equation}\label{ubMb}
  M_b(d) \le \left\{\begin{array}{ll}\frac{1}{3}d^2 + \frac{1}{3}d, & d \
      \textup{odd},\\[0.3em] 
      \frac{1}{3}d^2+\frac{1}{6}d, & d \ \textup{even}.\end{array}\right.
\end{equation}
$d$ of these cells are unbounded. 
This is a purely combinatorial result which is strongly related to the problem
of determining the maximum number of triangles in such arrangements which has
a long and rich history (see \cite{arrGoodman}). 
Notice that this bound is better than the one obtained by Kharlamov using
Hodge theory \cite{kharlOverview}.
It is known that the bound (\ref{ubMb}) is exact for infinitely many
values of $d$. 
The real folding polynomials $F_{\dR,d}^{A_2}(x,y)$ almost achieve this bound.  
Moreover, our arrangements have the very special property that all critical
points with a negative (resp.\ positive) critical value have the same critical
value $-1$ (resp.\ $+8$).  

To translate the upper bound on the number of black cells into an upper bound
on critical points we use the following lemma: 
\begin{Lem}[see Lemme 10, 11 in \cite{ortizAspects}]
  Let $f$ be a real simple line arrangement consisting of $d\ge3$ lines. 
  $f$ has exactly $d-1\choose2$ bounded open cells each of which contains
  exactly one critical point.  
  All the critical points of $f$ are non-degenerate.
\end{Lem}
It is easy to prove the lemma, e.g.\ by counting the number of bounded cells and
by observing that each such cell contains at least one critical point. 
Comparing this with the number $(d-1)^2-{d\choose2}={d-1\choose2}$ of all
non-zero critical points gives the result. 
Now we can show that our real line arrangements are asymptotically the best
possible ones for constructing surfaces with many singularities: 
\begin{Thm}\label{thmubcp}
  The maximum number of critical points with the same non-zero critical value 
  $0\ne v\in\dR$ of a real simple line arrangement is bounded by $M_b(d)-d$,
  where $d$ is the number of lines.  
  In particular, the maximum number of critical points on two levels of such
  an arrangement does not exceed ${d\choose2}+M_b(d)-d\approx \frac56 d^2$. 
\end{Thm}
\begin{proof}
  In view of the upper bound (\ref{ubMb}) for the maximum number $M_b(d)$ of
  black cells of a real simple line arrangement we only have to verify that
  any bounded cell contains only one critical point. 
  But this follows from the preceding lemma. 
\end{proof}

Chmutov showed a much more general result (\cite{chmuDefsCritPts}, see
\cite{chmuCritVals} for the case of non-degenerate critical points): 
For a plane curve of degree $d$ the maximum number of critical points on two levels
does not exceed  $\approx \frac{7}{8}d^2$. 
In \cite{chmuCritVals}, he conjectured $\approx \frac{5}{6}d^2$ to be the
actual maximum which is attained by the complex line arrangements
$F_{d}^{A_2}(x,y)$ he used for his 
construction (and also by the real line arrangement $F_{\dR, d}^{A_2}(x,y)$).  
Thus, our theorem \ref{thmubcp} is the verification of Chmutov's conjecture in the
particular case of real simple line arrangements.  
As Chmutov remarked in \cite{chmuP3}, such an upper bound immediately implies
an upper bound on the maximum number of nodes on a surface in separated
variables: 
\begin{Cor}
  A surface of the form $p(x,y) + q(z)=0$ cannot have more than $\approx
  \frac{1}{2}d^2\smcdot \frac12 d + \frac{1}{3}d^2\smcdot \frac12 d =
  \frac{5}{12}d^3$ nodes if $p(x,y)$ is a real simple line
  arrangement.  
  This number is attained by the surfaces
  $\textup{Chm}_{\dR,d}^{A_2}(x,y,z)$ defined in theorem \ref{thmChmRd}. 
\end{Cor}

Comparing this number to the upper bound $\approx \frac{5}{12}d^3$
on the zero$^\textup{th}$ Betti number (see e.g., \cite[p.\
533]{kharTopRealMan}) one is tempted to ask if it 
is possible to deform our singular surfaces to get examples with many real
connected components. 
But our surfaces $\textup{Chm}_{\dR,d}^{A_2}(x,y,z)$ only contain $A_1^-$
singularities which locally look like a cone ($x^2+y^2-z^2=0$). 
When removing the singularities from the zero-set of the surface every
connected component contains at least three of the singularities.
Thus, the zero$^\textup{th}$ Betti number of a small deformation of our surfaces are not
larger than $\approx \frac{5}{3\smcdot12}d^3$ which is far below the number
$\approx \frac{13}{36}d^3$ resulting from Bihan's construction
\cite{bihanBetti}.    

Conversely, we may ask if it is always possible to move the lines of a simple
real line arrangement in such a way that all critical points which have a
critical value of the same sign can be chosen to have the same critical
value. 
If this were true then it would be possible to improve our lower bound for the
maximum number $\mu^\dR(d)$ of real nodes on a real surface of degree $d$
slightly because it is known that the upper bounds for the maximum number
$M_b(d)$ of black cells are in fact exact for infinitely many $d$. 
E.g., in the already cited article \cite{harborthTwoCol}, Harborth gave an
explicit arrangement of $13$ straight lines which has $\frac{1}{3}\smcdot13^2 +
\frac{1}{3}\smcdot13 - 13=47$ bounded black regions.  
When regarding this arrangement as a polynomial of degree $d=13$ it has
exactly one critical point with a negative critical value within each of the 
black regions. 
If all these negative critical values could be chosen to be the same then such a
polynomial would lead to a surface with
${13\choose2}\smcdot\lceil\frac{13-1}{2}\rceil + 
47\smcdot\lfloor\frac{13-1}{2}\rfloor=750>732$ nodes.
Similarly, such a surface of degree $9$ would have $228>216$ nodes. 
In the case of degree $7$ the construction would only yield $96$ nodes which
is less than the number $99$ found in \cite{labs99}.

Notice that it is not clear that line arrangements are the best plane curves
for our purpose, and we may ask: 
Is it possible to exceed the number of critical points on two levels of the
line arrangements $F_{\dR,d}^{A_2}(x,y)$ using irreducible curves of higher
degrees?    
Either in the real or in the complex case?
This is not true for the real folding polynomials. 
E.g., those associated to the root system $B_2$ consist of many ellipses and
yield surfaces with fewer singularities (see 
\cite{sonjaChmuVars}).  

%%%%%%%%%%%%%%%%%%%%%%%%%%%%%%%%%
%
% the bibliography:
%
\nocite{labsAlgSurf}

\bibliographystyle{plain} 
\bibliography{papers}

\end{document}